\documentclass[11pt]{amsart}
\usepackage{amssymb,hyperref}

\newtheorem{bigthm}{Theorem}
\newtheorem{bigthmp}{Theorem}

\newtheorem{thm}{Theorem}[section]
\newtheorem{lem}[thm]{Lemma}
\newtheorem{prop}[thm]{Proposition}

\theoremstyle{definition}
\newtheorem{defn}[thm]{Definition}
\newtheorem{exmp}[thm]{Example}

\theoremstyle{remark}

\newtheorem*{rem*}{Remark}
\newtheorem{rem}[thm]{Remark}
\newtheorem*{rems*}{Remarks}

\numberwithin{equation}{section}

\def\bR{{\mathbb R}}  
 \def\bQ{{\mathbb Q}} \def\bZ{{\mathbb Z}}

\def\ens{\ensuremath} %%%%%%%%% ENSuremath
\newcommand\bydef{\ens{\stackrel{\text{def}}{=}}}% \bydef

\def\GG{\Gamma} 

\newcommand\hsp[1]{\mbox{}\hspace{#1mm}} % \hsp
\begin{document}
\bibliographystyle{amsplain}

\title{Almost everywhere convergence and polynomials}

\author{Michael Boshernitzan}

\address{Deparment of Mathematics, Rice University, Houston, TX~77005,
  USA}

\email{michael@rice.edu}

\author{M\'at\'e Wierdl}

\address{Department of Mathematical Sciences, University of Memphis,
  373 Dunn Hall, Memphis, TN 38152-3240}

\email{mw@csi.hu}

\date{\today}

\begin{abstract}
  Denote by $\GG$ the set of pointwise good sequences.  Those are
  sequences of real numbers $(a_k)$ such that for any measure
  preserving flow $(U_t)_{t\in\bR}$ on a probability space and for any
  $f\in L^\infty$, the averages $\frac{1}{n} \sum_{k=1}^{n}
  f(U_{a_k}x) $ converge almost everywhere.

  We prove the following two results.
  \begin{enumerate}
  \item If $f\!: (0,\infty)\to\bR$ is continuous and if
    $\big(f(ku+v)\big)_{k\geq 1}\!\in\GG$ for all $u, v>0$, then $f$
    is a polynomial on some subinterval $J\subset (0,\infty)$ of
    positive length.
  \item If $f\!: [0,\infty)\to\bR$ is real analytic and if
    $\big(f(ku)\big)_{k\geq 1}\!\in\GG$ for all $u>0$, then $f$ is a
    polynomial on the whole domain $[0,\infty)$.
  \end{enumerate}

  These results can be viewed as converses of Bourgain's polynomial
  ergodic theorem which claims that every polynomial sequence lies in
  \,$\GG$.

\end{abstract}

\maketitle

\tableofcontents
\section{Introduction}

For $1\leq p\leq\infty$, \ a sequence of real numbers $(a_k)_{k\geq1}$
is said to be $p$-good pointwise if for any measure preserving flow
$(U_t)_{t\in\bR}$ on a probability space $(X, {\mathcal B}, \mu)$ and
for any $f\in L^p(X,\mu)$, the averages $\frac{1}{n} \sum_{k=1}^{n}
f(U_{a_k}x) $ converge almost everywhere.

Denote by $\GG_p$ the set of $p$-good pointwise sequences and set
$\GG=\GG_\infty$.

J.~Bourgain proved in a series of papers (see \cite{bourgain} for the
proof and further references) that polynomial sequences lie in
$\GG_p$, for all $p>1$. Bourgain, in fact, formulated his result for
$\mathbb Z^d$ actions instead of $\mathbb R$ actions.  To see the idea how we
can translate the $\mathbb Z^d$ results to results on a flow, consider
the example $f(x)=\sqrt 2x^2 +x$.  The transformation
$T_{u,v}=U_{\sqrt 2u+v}$ is a $\mathbb Z^2$ action, and
$U_{f(n)}=T_{n^2,n}$, hence Bourgain's result applied to $T$ gives the
result for the flow $U_t$.  

For a simpler proof (in the case $p=2$) of Bourgain's result, and
for references see \cite{qw:squares}.  It is a famous open problem
whether or not $(k^2)\in\GG_1$. 

In the paper we prove the following two theorems.  (Each can be viewed
as a converse of J.~Bourgain's result mentioned above).

Denote \ $\bR^+=\{r\in\bR\mid r>0\}=(0,\infty)$, \ and \
$\overline{\bR}^+=\{r\in\bR\mid r\geq0\}=[0,\infty)$.

\begin{bigthm} \label{thm:1} Let $f\!:\overline{\bR}^+\to\bR$ be a
  real analytic function such that the set
  \begin{equation} \label{eq:U} U=\{u>0: \big(f(ku)\big)_{k\geq 1} \in
    \GG\}, \quad \text{ where } \quad \GG\bydef\GG_\infty,
  \end{equation}
  is uncountable. Then $f(x)$ must be a polynomial.
\end{bigthm}

\begin{bigthm} \label{thm:2} Let $U, V\subset\bR^+$ be two non-empty
  open subsets.  Let $f\!: \bR^+\to\bR$ be a real continuous function
  such that $\big(f(ku+v)\big)_{k\geq1}\!\in\GG$, for all $u\in U,
  v\in V$.  Then there exists a subinterval $J\subset \bR^+$ (of
  positive length) such that $f$ restricted to $J$ is a polynomial.
\end{bigthm}

We remark that, in general, the conclusion of Theorem~\ref{thm:2}
cannot be that $f(x)$ agrees with a polynomial on a halfline, say.
Indeed, we have

\begin{prop}\label{prop:vlast}
  The function
  $$
  f(x)=\mathbf{dist}(x,\bZ)=\min(x-[x], [x]+1-x), \quad x\in\bR,
  $$   
  satisfies the conditions of Theorem~\ref{thm:1} with $U=V=\bR^+$.
\end{prop}

The above proposition shows that $f(x)$ satisfying the conditions of
Theorem~\ref{thm:2} \ (even with \text{$U=V=\bR^+$})\, may be a
periodic non-constant function, so definitely not a polynomial on large
subintervals of \ $\bR^+$.

The claim of Proposition~\ref{prop:vlast} follows from the results in
\cite{bjw}.

Note that our proof of Theorem~\ref{thm:1} depends on the fact that
$f$ is analytic at $0$.  (See Remark~\ref{rem:1} in the next section).

The following result (first proved in \cite{bbb}) plays a central role
in the proofs of both theorems.

\begin{lem}[\cite{bbb}]\label{lem:1}
  Any sequence $(a_k)\in \GG$ must be linearly dependent over the
  field $\bQ$ of rational constants.
\end{lem}

For a simpler proof of this lemma (avoiding the use of Bourgain's
entropy method) we refer the reader to~\cite{jw}.

\section{Proof of Theorem~\ref{thm:1}}

\begin{proof}[Proof of Theorem~\ref{thm:1}] 
  In view of Lemma~\ref{lem:1}, for every $x\in U$, there are an integer
$n(x)\geq 1$ and an $n(x)$-tuple of rationals
$${\bf q}(x)=(q_1(x), \ldots  , q_{n(x)}(x) ) \in \bQ^{n(x)}, \quad q_{n(x)}=1,$$
such that
$$\sum_{k=1}^{n(x)} q_k(x) f(kx) = 0.$$

Since the set of possible pairs $\big(n(x),{\bf q}(x)\big)$ is
countable while the set of $x\in U$ is uncountable, there is a pair
$\big( n, {\bf q}\big)$ corresponding to an uncountable set $U'\subset
U$ of $x$:
$$\big(n(x), {\bf q}(x)\big) = 
\big(n, {\bf q}\big), \qquad \text{for } \ x\in U'.$$

Since $f(x)$ is analytic, the identity
\begin{equation}
  \label{eq:recurrence} 
  \sum_{k=1}^n q_k f(kx) = 0, \quad q_n=1,
\end{equation}
extends (from $x\in U'$) to all $x\in [0,\infty)$ (the set $U'$ being
uncountable has an accumulation point).

Set $f(x)=\sum_{\,r\geq 0} c_r x^r$ to be the series expansion of $f$
at $0$.  Then from the identity
$$  0=\sum_{k=1}^n q_k f(kx)=
\sum_{\,r\geq 0}\big(\sum_{k=1}^n q_k k^r \big)\, c_r x^r
$$
one deduces that, for all integers $r\geq 0$, either $c_r=0$ \ or \ $
\sum_{k=1}^n q_k k^r=0$. It follows that

$$  \sum_{k=1}^n q_k k^r=0,$$
for all $r$ in the set $K=\{r\in \bZ^+: c_r\neq 0\}. $

We observe that
$$ \lim_{r\to +\infty} \,\sum_{k=1}^n q_k k^r=+\infty$$
(the last term $q_n n^r=n^r$ is dominant in the sum).

We conclude that the set $K=\{r\in \bZ^+: c_r\neq 0\}$ is finite, and
$f(x)$ is a polynomial.
\end{proof}

\begin{rem}\label{rem:1}
  We don't know whether Theorem~\ref{thm:1} holds if the domain of $f$
  is assumed to be $(0,\infty)$ rather than $[0,\infty)$.  The
  recurrence relation \eqref{eq:recurrence} still holds in this
  setting. The conclusion of Theorem can be derived under the
  assumptions that $f$ is analytic on $(0,\infty)$ and that either $0$
  or $\infty$ is an isolated singularity of the analytic extension of
  $f$.
\end{rem}

\section{Proof of Theorem~\ref{thm:2}}

\begin{proof}[Proof of Theorem~\ref{thm:2}]
  In view of Lemma~\ref{lem:1}, for every $u\in U$ and $v\in V$ there
  is an integer $n(u,v)\geq1$ and an $n(u,v)$-tuple of rationals
$$
{\bf q}(u,v)=(q_1(u,v), q_2(u,v), \ldots , q_{n(u,v)}(u,v)) \in
\bQ^{n(u,v)}, \quad {\bf q}(u,v)\neq{\bf 0},
$$
such that
$$
\sum_{k=1}^n q_k(u,v)\, f(ku+v) = 0.
$$

For every integer \ $n\geq1$ \ and an $n$-tuple of rationals \ ${\bf
  q}\in\bQ^n$, denote
$$
K(n,{\bf q})=\big\{(u,v)\in U\times V \ \big| \ n(u,v)=n \hsp2 \text{
  {\Small and} } \hsp2 {\bf q}(u,v)={\bf q}\big\}.
$$

Since $U\times V$ is a countable union of its closed subsets $K(n,{\bf
  q})$, by the Baire category theorem there is a choice of \ $n$ \ and
\ ${\bf q}=(q_1, q_2, \ldots, q_n)\in\bQ^n$, with not all \ $q_k=0$,
such that the set $K(n,{\bf q})$ contains a non-empty interior, say
the set \ $U'\times V'$ \ where $U'\subset U$ and $V'\subset V$ are
non-empty open subintervals of $(0,\infty)$.  We conclude that

\begin{equation}\label{eq:kuv}
  \sum_{k=1}^n \, q_k\, f(ku+v) = 0, \quad \text{ for }  u\in U',\  v\in V'.
\end{equation}

\begin{defn}\label{def:quin}
  Let $U\subset\bR$ be an open set, and let $X\subset U$ be a finite
  subset, {\bf card}$(X)\geq1$.  Let $f: U\to\bR$, $g: X\to\bR$ and
  $h: X\to\bR$ be three functions such that $f$ is continuous and $g$
  is injective.  The quintuple $(f,g,h,X,U)$ is called balanced if
  \begin{equation}\label{eq:bal}
    \sum_{x\in X}  h(x) f(x+s+t\,g(x)) = 0,
  \end{equation}
  provided that $|t|$ and $|s|$ are small enough.
\end{defn}

\begin{exmp}
  Let $U', V', n\geq1$ and ${\bf q}=(q_1, q_2, \ldots, q_n)\in\bQ^n$
  be such as described in the paragraph preceding \eqref{eq:kuv}.
  Pick $u_0\in U', v_0\in V'$ and let $X$ be the set $X=\{x_k \mid
  1\leq k\leq n\}$ where $x_k=ku_0+v_0$.  Define $g, h: X\to \bR$ as
  follows: $g(x_k)=k, \, h(x_k)=q_k$.  With these choices, it follows
  from \eqref{eq:kuv} that the quintuple $(f,g,h,X,\bR^+)$ is
  balanced.

  Indeed, by setting $u=u_0+t$ and $v=v_0+s$, we obtain
  $$
  \sum_{x\in X} h(x) f(x+s+t\,g(x))=\sum_{k=1}^n h(x_k)
  f(x_k+s+t\,g(x_k))=\sum_{k=1}^n q_k\,f(x_k+s+tk)=
  $$

  $$
  =\sum_{k=1}^n q_k \,f(x_k+s+tk)=\sum_{k=1}^n q_k
  \,f\big((v_0+s)+(u_0+t)k\big)= \sum_{k=1}^n q_k \,f(v+uk)=0,
  $$
  since $u\in U',\, v\in V'$ if both $t, s$ are close to $0$.
\end{exmp}

Now the claim of Theorem~\ref{thm:2} follows from the following
result.

\begin{prop}\label{prop:last}
  Let\, $(f,g,h,X,U)$ be a balanced quintuple in the sense of
  Definition~\ref{def:quin}, with {\bf card}$(X)=n\geq1$.  Let\,
  $x_0\in X$ be such that\, $h(x_0)\neq0$.  Then $f(x)$ is a
  polynomial of degree $\leq n-2$ in a neighborhood of\, $x_0$.
\end{prop}

By definition, the zero constant is a polynomial of degree $-1$.
\end{proof}

\section{Proof of Proposition~\ref{prop:last}; smooth case}

\begin{proof}[Proof of Proposition~\ref{prop:last}; smooth case]
  First we provide proof under the additional assumption that $f\in
  C^\infty(U)$.  Without loss of generality, we may assume that \
  $0\notin h(X)$.  The proof is by induction on \ $n=\,${\rm
    card}$(X)$.

  If $n=1$, then by setting $t=0$ we get $f(x_0+s)=0$, for all $|s|$
  small enough, i.\,e.  $f(x)$ vanishes in a neighborhood of
  $x_0$. The case $n=1$ is validated.

  Next we assume that $n\geq1$, that $X=\{x_0, x_1, \ldots, x_n\}$ and
  that the claim of proposition is validated if card$(X)\leq n$.
  Observe that for an arbitrary real constant \ $c$ \ a quintuple
  $(f,g,h,X,U)$ is balanced if and only if a quintuple $(f,g+c,h,X,U)$
  is.  This is because $ f(x+s+t\,g(x))=f(x+s-ct+(g(x)+c)t), $ and the
  pair $(s,t)$ is close to $(0,0)$ (in the $\bR^2$ metric) if and only
  if $(s-ct, t)$ is.

  Thus we may assume that $g(x_n)=0$ (after replacing $g(x)$ by
  $g(x)-g(x_n)$).  Thus
$$
0=\sum_{x\in X} h(x) f(x+s+t\,g(x))=\sum_{k=0}^{n-1} h(x_k)
f(x_k+s+t\,g(x_k))+h(x_n)f(x_n+s).
$$

Taking the partial derivative $\frac{\partial}{\partial t}$ we get
$$
\sum_{k=0}^{n-1} h(x_k)g(x_k) f'(x_k+s+t\,g(x_k))=0
$$
and conclude that the quintuple \ $(f',g|_{X'},hg|_{X'},X',U)$ \ is
balanced where \ $X'=\{x_0,\ldots, x_{n-1}\}=$ $=X\setminus \{x_n\}$.
Moreover, $g|_{X'}$ is injective, and $0\notin hg(X')$.

By the induction hypothesis, for any $x'\in X'$, the derivative $f'$
is a polynomial of degree $\leq n-2$ in some neighborhood of $x'$.  It
follows that $f(x)$ is a polynomial of degree $\leq n-1=(n+1)-2$ in a
neighborhood of $x_0$.  This completes the proof of Proposition
\ref{prop:last} under the added assumption that $f\in
C^\infty(\bR^+)$.
\end{proof}

\section{Proof of Proposition~\ref{prop:last}; continuous case}

\begin{proof}[Proof of Proposition~\ref{prop:last}; continuous case]
  Since $(f,g,h,X,U)$ is a balanced quintuple, there exists an
  $\epsilon>0$ such that \eqref{eq:bal} holds provided that \ $|s|,
  |t|<\epsilon$.  In the preceding section we proved that, under the
  additional condition that $f\in C^\infty(\bR^+)$, there exists a
  neighborhood $W$ of a point $x_0$ such that $f|_W$ is a polynomial
  of degree $\leq n-2$.

  Our proof provides slightly more: This neighborhood $W$ depends only
  on $\epsilon, g, h, X$ and $U$ but not on $f\in C^\infty(\bR^+)$.
  The observation will be used in what follows.

  Now we move to the general case of $f\in C(\bR^+)$ (rather than
  $f\in C^\infty(\bR^+)$).

  Fix any function $\phi\in C^\infty(\bR)$ which satisfies\, $\int_\bR
  \, \phi(t)\, dt=1$ and vanishes outside the interval $[-1,1]$.  For
  integers $m\geq1$, denote
$$
\phi_m(x)=m\phi(mx), \qquad f_m(x)\bydef\big(\phi_m\circ
f\big)(x)=\int_\bR \phi_m(t)f(x-t)\, dt.
$$

The sequence of kernels $\phi_m(x)$ is known to converge the
$\delta$-function at $0$ in the sense that
\begin{equation}\label{eq:conv}
  \lim_{m\to\infty} f_m(x)=f(x), \quad \text{ for all } \  x\in U.
\end{equation}  
Note that all $f_m(x)\in C^\infty(U_m)$ where
$$
U_m\bydef\big\{x\in U \mid \text{\rm dist}(x,\partial
U)>\tfrac1m\big\},
$$
and\, dist$(x,\partial U)$\, stands for the distance between $x$ and
the boundary of\, $U\subset \bR$.

Let $V$ be a neighborhood of the set $X$ such that its closure
$\overline V$ is compact and is contained in $U$.  Then the pointwise
convergence in \eqref{eq:conv} is uniform on $V$, and in fact $f_m\in
C^\infty(V)$ holds for all $m$ large enough.  It is also clear that
for $m$ large enough, $(f_m,g,h,X,V)$ forms a balanced quintuple
because $(f,g,h,X,V)$ does; moreover, there exists $\epsilon'>0$ such
that
 $$
 \sum_{x\in X} h(x)\,f_m(x+s+tg(x))=0
 $$
 holds simultaneously for all large $m$ (say, $m>m_0$) and
 $s,t\in(-\epsilon',\epsilon')$.  It follows that there exists a
 neighborhood $W\subset V$ of a point $x_0$ such that each function
 $f_m$ is a polynomial of degree $\leq n-2$ in it.  (Here we use the
 observation made in the second paragraph of this section).  In view
 of the uniform convergence \eqref{eq:conv}, $f|_W$ is also a
 polynomial of degree $\leq n-2$, completing the proof of
 Proposition~\ref{prop:last}.
\end{proof}

 \section{Concluding remarks}
 
 The following is a slightly more general version of
 Theorem~\ref{thm:1}.

 \begin{bigthmp} \label{thm:3} Let $(r_k)_{k\geq1}$ be a sequence of
   distinct positive numbers, let $f\!:\overline{\bR}^+\to\bR$ be a
   real analytic function such that the set
   \begin{equation} \label{eq:U2} U=\{u>0: \big(f(r_ku)\big)_{k\geq 1}
     \in \GG\}, \quad \text{ where } \quad \GG\bydef\GG_\infty,
   \end{equation}
   is uncountable. Then $f(x)$ must be a polynomial.
 \end{bigthmp}

 The proof of Theorem \ref{thm:3} is very similar to the proof of
 Theorem \ref{thm:1}.

\bibliography{polypoint}

\end{document}